\documentclass[12pt]{article}
\begin{document}

\begin{center}
\textbf{Generalized Nambu dynamics and vectorial Hamiltonians}
\bigskip

\textbf{V.N.Dumachev}\footnote{Voronezh Institute of the MVD of the Russia} \\
dumv@comch.ru
\end{center}

\begin{small}

On the basis of Liouville theorem the generalization of the Nambu
mechanics is considered.  Is shown, that Poisson manifolds of
n-dimensional multi-symplectic phase space have inducting by (n-1)
Hamiltonian k-vector fields, each of which requires of
(k)-hamiltonians.

\end{small}

\bigskip
\textbf{1.} Consider a system of the differential equations as
submanifold $\Sigma$ in a jet bundles $J^n(\pi)$: $E \to M $, and
[1]
\[
F(t,x_0,x_1,...,x_n)=0,
\]
where $t \in M \subset R$, $u = x_0 \in U \subset R$, $x_i \in J^i(
\pi) \subset R^n$, $E=M \times U$ . We select from $\Sigma \subset J
^n(\pi)$ the equations with Poisson structure and bracket
\begin{equation}
\label{eq1}\left\{{H,G} \right\}=X_H\rfloor dG=\mathcal{L}_{X_
{H}}G,
\end{equation}
and write Cartan distribution as
\begin{equation}
\label{eq2} \theta _i=dx_i-\{H,x_i\}dt.
\end{equation}
Here $\mathcal{L} _ {X _{H}} $ is the Lie derivative along vector
field $X _ {H} \in \Lambda ^ {1} $, $ \Lambda^n$ a exterior graded
algebra of $k$-vector fields, $H=H(x)$ - while unknown function.
Contact vector field ($X_H \rfloor \theta =0$) has the form
\begin{equation}
\label{eq3} X_H=\frac{\partial}{\partial t}+\{H,x_i\}\frac{\partial
}{\partial x_i}.
\end{equation}

\textbf{2.} Consider classical symplectic mechanics on $J^1(\pi)$.
Vector field $X_H^1$ on symplectic manifold $(M, \omega) $ is called
Hamiltonian, if the 1-form $\Theta =X_H^1\rfloor \Omega $ is closed
$d\Theta = 0 $ and exact (for contractible manifolds) [2]. It allows
to find a Hamiltonian $H$. In this case symplectic form is
\[
\Omega = dx _ 0 \wedge dx _1
\]
and
\[
\Theta = X _ {H} ^ {1} \rfloor \Omega = dH,
\]
where
\[
H=\frac{1}{2}\left( x_0^2+x_1^2 \right).
\]
Poisson structure (\ref{eq1}):
\[
\{F,G\}=\frac{\partial F}{\partial x_0} \cdot \frac{\partial
G}{\partial x_1}-\frac{\partial F}{\partial x_1} \cdot
\frac{\partial G}{\partial x_0}.
\]
a inducted by Hamiltonian vector field:
\[
X_H^1=\frac{\partial H}{\partial x_1} \cdot \frac{\partial
}{\partial x_0}-\frac{\partial H}{\partial x_0} \cdot
\frac{\partial}{\partial x_1}.
\]
Cartan distribution defined dynamic Hamiltonian equations in
canonical coordinates
\[
\frac{\partial x_i}{\partial t}=\{H, x_i\}.
\]
Volume of the Cartan differential forms (\ref{eq2})
\[
I = \theta _0 \wedge \theta _1=\Omega - X^1\rfloor \Omega \wedge dt
\]
give as Poincar\'e integral invariant $i \quad \left( {I = di}
\right)$:
\[
i = \frac{1}{2}(x_0 dx_1-x_1 dx_0)-H \wedge dt.
\]

According to the Liouville theorem anyone Hamiltonian field
conserved of the Volume form, i.e. Lie derivative of the 2-form
along vector field $X^1_H$ is zero: $\mathcal{L}_X\Omega=0$. In
other words, generated by a vector fields $X^1_H$ the
one-parametrical group symplectic transformations $\{g_t\}$ (phase
flow) preserves the canonical two-form $\Omega$, i.e. $g_t\Omega=0$.

\bigskip
\textbf{3.} Now, we have extended the previous calculations on
$J^2(\pi)$. From Liouville theorem we have 3-form
\[
\Omega=dx_0 \wedge dx_1 \wedge dx_2
\]
volume of the phase space.

\textbf{Theorem 1.} The volume 3-form  $ \Omega \in \Lambda ^3$ are
supposed by existence of 2 polyvector Hamiltonian fields: $X^1\in
\Lambda ^1$ and $X^2\in \Lambda ^2$.

\textbf{Proof:} On definition
\[
\mathcal{L}_X \Omega = X \rfloor d\Omega + d(X\rfloor \Omega)=0.
\]
Since $\Omega \in \Lambda ^ 3$, we see that $d\Omega = 0$, and
\[
d(X\rfloor \Omega)=0.
\]
From Poincar\'e lemma it follows that $X\rfloor \Omega $ is exact,
and
\[
X\rfloor \Omega = \Theta = dH.
\]

\noindent 1) If $\textbf {X}_H^1 \in \Lambda ^1$, then $ \Theta \in
\Lambda ^ 2 $, $H \in \Lambda ^ 1 $:
\[
X_H^1=\left( \frac{\partial H}{\partial x_1}-\frac{\partial
H}{\partial x_2} \right)\frac{\partial }{\partial x_0}+\left(
\frac{\partial H}{\partial x_2}-\frac{\partial H}{\partial x_0}
\right)\frac{\partial}{\partial x_1}+\left( \frac{\partial
H}{\partial x_0}-\frac{\partial H}{\partial x_1}
\right)\frac{\partial}{\partial x_2},
\]
and
\[
\Theta = X_H^1 \rfloor \Omega = \sum\limits_{i=0}^2 {dH \wedge
dx_i}.
\]
Poisson bracket $X_H^1 \rfloor dG = \{ H,G \}$ has the form
\[
\{ H,G \} = \left( \frac{\partial H}{\partial x_1}-\frac{\partial
H}{\partial x_2} \right)\frac{\partial G}{\partial x_0}+\left(
\frac{\partial H}{\partial x_2}-\frac{\partial H}{\partial x_0}
\right)\frac{\partial G}{\partial x_1}+\left( \frac{\partial
H}{\partial x_0}-\frac{\partial H}{\partial x_1}
\right)\frac{\partial G}{\partial x_2}
\]

\noindent 2) If $\textbf {X}_H^2 \in \Lambda ^2$, then  $ \Theta \in
\Lambda ^1 $, $H \in \Lambda ^0 $,
\[
X_H^2=\frac{1}{2}\left(\frac{\partial H}{\partial x_0} \cdot
\frac{\partial}{\partial x_1}\wedge \frac{\partial}{\partial x_2}+
\frac{\partial H}{\partial x_1} \cdot \frac{\partial}{\partial x_2}
\wedge \frac{\partial}{\partial x_0}+\frac{\partial H}{\partial x_2}
\cdot \frac{\partial}{\partial x_0} \wedge \frac{\partial} {\partial
x_1}\right),
\]
and
\[
\Theta = X_H^2 \rfloor \Omega = dH.
\]
Poisson bracket $X_H^2 \rfloor (dF \wedge dG)=\{H,F,G\}$ has the
form
\begin{eqnarray*}
\{F,G,H\}&=&\frac{1}{2}\left[ \frac{\partial H}{\partial x_0} \cdot
\left( \frac{\partial F}{\partial x_1}\frac{\partial G}{\partial
x_2} - \frac{\partial F}{\partial x_2}\frac{\partial G}{\partial
x_1}\right) \right. \\
\\
&+&\quad \; \left.\frac{\partial H}{\partial x_1} \cdot \left(
\frac{\partial F}{\partial x_2}\frac{\partial G}{\partial
x_0}-\frac{\partial F}{\partial x_0}\frac{\partial G}{\partial x_2}
\right)+ \frac{\partial H}{\partial x_2} \cdot \left(\frac{\partial
F}{\partial x_0}\frac{\partial G}{\partial x_1}-\frac{\partial
F}{\partial x_1}\frac{\partial G}{\partial x_0} \right) \right]
\end{eqnarray*}

\bigskip
\textbf{4.} Let's consider of generalization of the previous
calculations on $J^n(\pi)$. Let n-form

\[
\Omega = dx _ {0} \wedge dx _ {1} \wedge ... \wedge dx _ {n - 1}
\]

be Liuville's volume  of a phase space.

\textbf{Theorem 2.} Conservation laws for the Liuville's volume $
\Omega \in \Lambda ^ {n} $ suppose existence of $n-1$ polyvector's
Hamilton's fields $X^k \in \Lambda ^k \quad (1 \le k \le n-1)$.

\textbf{Proof:} On definition

\[
L_X \Omega = X \rfloor d\Omega + d(X\rfloor \Omega )=0.
\]

Since $ \Omega \in \Lambda ^n$, we see that $d\Omega =0$, and
$d(X\rfloor \Omega)=0$. From Poincare's lemma it follows that
$X\rfloor \Omega $ - is exact, and

\[
X\rfloor \Omega = \Theta = dH.
\]

\noindent If $\textbf {X} _H^k \in \Lambda ^k$, then $\Theta \in
\Lambda ^ {n - k} $, $H \in \Lambda ^ {n - k - 1} $ and $(1 \le k
\le n-1)$:

\[
X_{H}^{n - 1} = \frac{{1}}{{\left( {n - 1}
\right)!}}\sum\limits_{k = 0}^{n - 1} {\frac{{\partial
H}}{{\partial x_{k}} } \cdot \frac{{\partial }}{{\partial x_{0}} }
\wedge \frac{{\partial} }{{\partial x_{1}} } \wedge \left[
{\frac{{\partial} }{{\partial x_{k}} }} \right] \wedge ... \wedge
\frac{{\partial} }{{\partial x_{n - 1}} }} ,
\]

\[
X_{H}^{n - 2} = \frac{{1}}{{\left( {n - 2}
\right)!}}\sum\limits_{i < k}^{n - 1} {\left( {\frac{{\partial
H}}{{\partial x_{i}} } - \frac{{\partial H}}{{\partial x_{k}} }}
\right) \cdot \frac{{\partial} }{{\partial x_{0}} } \wedge \left[
{\frac{{\partial} }{{\partial x_{i}} } \wedge \frac{{\partial
}}{{\partial x_{k}} }} \right] \wedge ... \wedge \frac{{\partial
}}{{\partial x_{n - 1}} }} ,
\]

\[
X_{H}^{n - 3} = \frac{{1}}{{\left( {n - 3}
\right)!}}\sum\limits_{i < k < l}^{n - 1} {\left( {\frac{{\partial
H}}{{\partial x_{i}} } - \frac{{\partial H}}{{\partial x_{k}} } +
\frac{{\partial H}}{{\partial x_{l}} }} \right) \cdot
 \left[ {\frac{{\partial }}{{\partial
x_{i}} } \wedge \frac{{\partial} }{{\partial x_{k}} } \wedge
\frac{{\partial} }{{\partial x_{l}} }} \right] \wedge ... \wedge
\frac{{\partial} }{{\partial x_{n - 1}} }} ,
\]

…

\noindent for which

\[
\Theta ^{1} = X_{H}^{n - 1} \rfloor \Omega = dH,
\]

\[
\Theta ^{2} = X_{H}^{n - 2} \rfloor \Omega = \sum\limits_{i}^{n -
1} {dH \wedge dx_{i}} ,
\]

\[
\Theta ^{3} = X_{H}^{n - 3} \rfloor \Omega = \sum\limits_{i <
j}^{n - 1} {dH \wedge dx_{i} \wedge dx_{j}}  ,
\]

\[
\Theta ^{4} = X_{H}^{n - 4} \rfloor \Omega = \sum\limits_{i < j <
k}^{n - 1} {dH \wedge dx_{i} \wedge dx_{j} \wedge dx_{k}}  ,
\]

…

\[
\Theta ^{n - 1} = X_{H}^{1} \rfloor \Omega = \sum\limits_{i < j <
... < k}^{n - 1} {dH \wedge dx_{i} \wedge dx_{j} \wedge ... \wedge
dx_{k}}  .
\]

\textbf{Proposition.} For any $n$-form $\Omega $ and Hamiltonian
polyvector field $X_H^k \in \Lambda ^k$ existence Poisson Structure
$X_H^1\rfloor dx_i=\{H,x_i \}$, containing $k$ Hamiltonians:
\[
X_H^k \rfloor ( {dF}_{1}{\wedge dF}_{2}{\wedge }...{\wedge dF}_{k})
=\{H, F_1,F_2,...,F_k\}.
\]

\bigskip
\textbf{5.} Example. Euler's equations for rigid body in Poisson's
form may also be written as
\[
\left\{
\begin{array}{l}
\dot x=\;\;\; y-z, \\
\dot y=-x+z,\\
\dot z=\;\;\;x-y,
\end{array} \right.
\]
The vectorial form of this equation we write as
\[
\dot \textbf{x}=\textbf{Dx},
\]
where
\begin{center}
$\textbf{x}= \left(\begin{array}{l} x \\ y\\ z \end{array} \right)
\quad$ and $\quad \textbf{D}=\left(
\begin{array}{rrr}
0&1&-1 \\
-1&0&1\\
1&-1&0
\end{array} \right)$.
\end{center}
This vector flow is called Hamilton's flow if
\[
div \textbf{Dx}=0.
\]
This implies, that
\[
\textbf{Dx}=rot \textbf{h}.
\]
These expressions are a requirement of a closure of differential
form $d \omega=0\;$:
\[
\omega =(y-x)dy \wedge dz+(-x+z)dz \wedge dx+(x-y)dx \wedge dy=(rot
\textbf{h} \cdot d \textbf{S}).
\]
Using an Homotopy operator we have $\omega=d \nu\;$:
\begin{center}
$\nu =(\textbf{h} \cdot d \textbf{x}) \quad$ where $\quad
\textbf{h}=\frac{1}{3}\left(
\begin{array}{l}
y^2+z^2-x(y+z) \\
z^2+x^2-y(z+x)\\
x^2+y^2-z(x+y)
\end{array} \right),\quad
d\textbf{x}=\left( \begin{array}{l} dx \\ dy\\ dz \end{array}
\right)$,
\end{center}
or
\[
\nu =(y^2+z^2-x(y+z))dx+(z^2+x^2-y(z+x))dy+(x^2+y^2-z(x+y))dz.
\]
Hamilton's flow has the form
\[
X_h^1=\left( rot \textbf{h} \cdot \frac{\partial}{\partial
\textbf{x}}\right)=\left( \frac{\partial h_3}{\partial
y}-\frac{\partial h_2}{\partial z}
 \right)\frac{\partial}{\partial
x}+ \left( \frac{\partial h_1}{\partial z}-\frac{\partial
h_3}{\partial x}
 \right)\frac{\partial}{\partial
y} + \left( \frac{\partial h_2}{\partial x}-\frac{\partial
h_1}{\partial y}
 \right)\frac{\partial}{\partial
z},
\]
and Poisson bracket  $\left\{{H,G} \right\}=X_{H}\rfloor dG$ gives:
\[
\dot \textbf{x}=\left\{{\textbf{h},\textbf{x}} \right\},
\]
where $\textbf{h}$ - vectorial Hamiltonian.

From Lax pair
\[
\dot L=[ML], \qquad L=\left(
\begin{array}{lll}
x&z&y \\
z&y&x\\
y&x&z
\end{array} \right), \qquad
M=\frac{1}{2}\left(
\begin{array}{rrr}
0&-1&1 \\
1&0&-1\\
-1&1&0
\end{array} \right)
\]
for these flow we get two scalar invariants
\[
I_1=trL=x+y+z, \qquad
I_2=\frac{1}{2}\;trL^2=\frac{3}{2}(x^2+y^2+z^2),
\]
gives bivectors Hamilton's flow
\[
X_I^{2} = \frac{{1}}{{2}}\left( {\frac{{\partial I}}{{\partial
x_{0}} } \cdot \frac{{\partial} }{{\partial x_{1}} } \wedge
\frac{{\partial }}{{\partial x_{2}} } + \frac{{\partial
I}}{{\partial x_{1}} } \cdot \frac{{\partial} }{{\partial x_{2}} }
\wedge \frac{{\partial} }{{\partial x_{0}} } + \frac{{\partial
I}}{{\partial x_{2}} } \cdot \frac{{\partial }}{{\partial x_{0}} }
\wedge \frac{{\partial} }{{\partial x_{1}} }} \right).
\]
For two Hamiltonians $I_1,I_2$ it gives $X_{I_1I_2}^1=X_{I_1}^2
\rfloor dI_2$:
\[
X_{I_1I_2}^1 = \left( \frac{\partial I_1}{\partial y}\frac{\partial
I_2}{\partial z}-\frac{\partial I_1}{\partial z}\frac{\partial
I_2}{\partial y}
 \right)\frac{\partial}{\partial
x}+ \left( \frac{\partial I_1}{\partial z}\frac{\partial
I_2}{\partial x}-\frac{\partial I_1}{\partial x}\frac{\partial
I_2}{\partial z}
 \right)\frac{\partial}{\partial
y} + \left( \frac{\partial I_1}{\partial x}\frac{\partial
I_2}{\partial y}-\frac{\partial I_1}{\partial y}\frac{\partial
I_2}{\partial x}
 \right)\frac{\partial}{\partial
z},
\]
and Poisson bracket
\[
X_{H}^{2} \rfloor \left( {dF \wedge dG} \right) = \left\{ {H,F,G}
\right\}
\]
with dynamics equations
\[
\dot \textbf{x}=\{\textbf{I}_1,\textbf{I}_2,\textbf{x}\}.
\]

 \bf{References} \rm

1.  Griffiths P., Exterior Differential Systems and the Calculus of
Variations, Birk\"auser, Boston, 1983.

2. Nambu Y. Generalized Hamiltonian dynamics // Phys.Rev.D, 1973,
V.7, N.8, P.5405-5412.

\end{document}